\documentstyle[12pt]{article}

\begin{document}

\newtheorem{theo}{Theorem}[section]
\newtheorem{lemm}{Lemma}[section]
\newtheorem{prop}{Proposition}[section]

\newcommand{\pr}{{\bf Proof. }} 
\newcommand{\rr}{\Bbb{R}}
\newcommand{\wh}{\widehat}
\newcommand{\wt}{\widetilde}
\newcommand{\ol}{\overline} 
\newcommand{\ra}{\rightarrow}
\newcommand{\la}{\lambda}
\newcommand{\vep}{\varepsilon}

\title{ON FLAT CONNECTIONS INDUCED OVER COVERING MAPS}
\author{Dionyssios Lappas}
\date{}
\maketitle

\begin{abstract}
\footnote{This work is in final form and no part of this will publish 
elsewhere.}
Flat connections induced over covering maps are studied and the trivial ones
among them are discribed. In the sequel, we deal with the resulting
holonomy bundles.
\end{abstract}

\ \\
{\bf (AMS) subject Classification:} 53C05, 57R22.

\setcounter{section}{-1}
\section{Introduction.}

Let $P,M$ be connected paracompact smooth manifolds and suppose that
$P(M,G)$ is a principal fiber bundle, where the Lie group $G$ acts smoothly
from the left on the total space $P$. Given the flat connection $\omega$ on
$P(M,G)$ and a regular covering manifold $(\wh{M},q,M)$ of $M$, a flat
connection $\wh{\omega}$ is constructed on the bundle $\wh{P}(\wh{M},G)$,
induced over $\wh{M}$. Furthermore, over the universal covering manifold
$\wt{M}$ the induced connection is trivial, i.e., it is (isomorphic to) the
flat connection on the product bundle $\wt{P}=G\times \wt{M}$, \cite{K2}.

In  this work {\it we are looking for all those coverings, where the induced
connection is trivial}. In this respect, we compute the holomony group of an
induced connection in terms of the holomony morphism defined by the given
connection and the fundamental group of the covering. To be more precise,
in Section 1 we prove the following.
\begin{theo}
Let $h_\omega :\pi_1 (M) \rightarrow G$ denote the holomony homorphism for
the flat connection $\omega$ of the principal bundle $P(M,G)$ and
$\pi_1(\wh{M})$ be the fundamental group of the covering $(\wh{M},q,M)$.
Then, for the induced by $q$ connection $\wh{\omega}$, we have the formula
$Im(h_{\wh{\omega}}) \cong (h_{\omega}|\pi_1(\wh{M}))$.
\end{theo}
\newtheorem{cor0}[theo]{Corollary}
\begin{cor0}
The induced connection is trivial if and only if $\pi_1(\wh{M}) \subseteq
Kerh_{\omega}$. 
\end{cor0}

In Section 2 we deal with certain coverings inducing non trivial
connections and study the corresponding {\it holomony bundles}. We prove
that if $Kerh_{\omega} \subseteq \pi_1(\wh{M})$, all these bundles have the
same total space, namely that of $P(M,G)$ (cf. 2.2).

\section{The holomony morphism for the induced bundles.}

We first formulate a preparatory result concerning bundles induced over
covering projections.
\begin{lemm}
Let $(E,p,B,F)$ be a (locally trivial) fiber bundle with connected
fiber, $(B_1,f,B)$ a covering space and $(E_1,p,B_1,F)$ denote
the bundle induced by $f$. Then the induced map $f_{\xi}:E_1 \rightarrow E$
is a covering map.
\end{lemm}
\pr Using standard notation (see \cite{H1}, Chapter 2), we have the following
commutative diagram.
$$
\begin{array}{ccccc}
&&f_{\xi}&&\\
&E_1&\longrightarrow & E &\\
p_1&\downarrow&&\downarrow&p\\
&B_1&\longrightarrow&B&\\
&&f&&
\end{array}
$$
Let $V$ be an evenly covered neighbourhood in $B$ which is trivializable
with respect to $p$, i.e., $f^{-1}(V)=\cup V_i$ and $p^{-1}(V)\cong
F\times V$, where the union is discrete. In view of the above diagram, we
have $f_{\xi}^{-1}(p^{-1}(V))=p^{-1}_1(\cup V_i)=\cup p_1^{-1}(V_i)$.
Therefore, the open set $W=p^{-1}(V)$ is the desired evenly covered
neighbourhood for $f_{\xi}$.

We now come on our main object of study and fix first some notation. Given
the principal bundle $P(M,G)=(P,p,M,G)$ and the covering $(\wh{M},q,M)$, it
follows from Lemma 1.1 that the total space $\wh{P}$ is a covering space of
$P$. Denoting by $\wh{p}$ and $\wh{q}$ the corresponding projections, we
have the following commutative diagram.
$$
\begin{array}{ccccc}
&&&\wh{p}&\\
G&\longrightarrow&\wh{P}&\longrightarrow&\wh{M}\\
&&\wh{q} \downarrow&&\downarrow q\\
G&\longrightarrow&P&\longrightarrow&M\\
&&&p&
\end{array}
$$
If $\omega$ is a flat connection on $P(M,G)$ and $\wh{\omega}$ its lifting on
$\wh{P}(\wh{M},G)$, then $\wh{\omega}$ is also flat and $\wh{q}$ {\it preserves
the horizontality defined by these connections}. Let $N(x_o)$ be the maximal
horizontal leaf of $\omega$ through a basic point $x_o\in P$ and
$\wh{x}_o\in \wh{q}^{-1}(x_o)$. Then $N(\wh{x}_o)$ {\it is a covering space 
of} $N(x_o)$, {\it the restriction} $p|N(x_o):N(x_o) \rightarrow G/P \cong
M$ {\it is a covering map and} $N(x_o)$ {\it is the total space for the
holomony bundle of} $\omega$. For further terminology see \cite{K2},
\cite{T3}; Chapter 4 and \cite{V4} where proofs of these facts and exact
references are given. In the sequel base points, if not necessary, would be
ommited.
\newtheorem{pro2}[lemm]{Proposition}
\begin{pro2}
Let $h_{\omega}$ (resp. $h_{\wh{\omega}}$) be the holomony morphism of
$\omega$ (resp. $\wh{\omega}$). Then, $h_{\wh{\omega}}=h_{\omega}\circ
q_{\#}$, where $q_{\#}=q_{\#,x_o}:\pi_1(\wh{M},\wh{p}(\wh{x}_o)) \rightarrow
\pi_1 (M,p(x_o))$ is the induced by $q$ morphism, at the fundamental group
level. 
\end{pro2}
\pr Let $w$ be a path in $\wh{M}$ closed at $\wh{p}(\wh{x}_o)$ and
$r=q\circ w$. If $\ol{w}$ denotes the horizontal lifting of $w$ starting at
$\wh{x}_o$, then there exists {\it exactly one} $g_w\in G$ such that
$\ol{w}(1)=h_{\wh{\omega}}([w])= g_w\wh{x}_o$. Using similar notation, we
have $\ol{r}(1)=h_{\omega}([r])=g_rx_o$. We claim that $g_r=g_w$. Indeed,
since $\wh{q}$ preserves horizontality, the path $\wh{q}\circ \ol{w}$ is
also the horizontal lifting of $r$ starting at $x_o$. Therefore,
$$\wh{q}(g_w\wh{x}_o)=(\wh{q}\circ \ol{w})(1)=\ol{r} (1)=g_rx_o.$$
Because $\wh{q}$ also preserves the orbit parametrization and the fibration
is principal, it follows that $g_wx_o=g_rx_o$ and finally $g_w=g_{r}$, 
as required.

\smallskip
\indent
As $q_{\#}$ is a monomorphism, the proof of the Theorem follows from Prop.
1.2, as well as the Corollary too. The further discription of
$\wh{P}(\wh{M},G)$ is carried out below.

\section{The holomony bundles for certain induced connections.}

We now proceed further and deal with the resulting holomony bundles over
certain coverings of $M$. Our previous notation is still in force.
\begin{prop}
If for the regular covering $(\wh{M},q,M)$ the relation
$Kerh_{\omega}$ $=\pi_1(\wh{M})$ holds, then the flat connections $\omega$ and
$\wh{\omega}$ have the same total space for their holonomy bundles, i.e.,
$N(x_0) \cong N(\wh{x}_o)$.
\end{prop}
\pr As $N(\wh{x}_o) \cong \wh{M}$ and it is a covering space of $N(x_o)$,
it is enough to prove that they have isomorphic fundamental groups. Let $z$
be a horizontal path closed at $x_o,\wh{z}$ its lift with respect to the
covering $(N(\wh{x}_o),\wh{q}|N(\wh{x}_o)$, $N(x_o))$ starting from $\wh{x}_o$
and suppose that $\wh{z}(1) \neq \wh{z}(0)$. Since $\wh{P}(\wh{M},G)$ is
the product bundle and has trivial holonomy, the path $\wh{p}\circ \wh{z}$
is not closed in $\wh{M}$, while $(q \circ \wh{p})\circ \wh{z}=p\circ z$ is
closed in $M$. Because $(\wh{M},q,M)$ is a covering, it follows that the
homotopy class $[p\circ z]$ is not contained in $Kerh_{\omega}\cong
\pi_1(\wh{M})$. This in particular implies that the horizontal lift $\ol{p \circ
z}$ of $p\circ z$, with respect to $p$, is not closed. As $z$ is also
horizontal, we have $\ol{p \circ z}=z$, which is closed, contradiction.
Thus, for every $[\gamma]\in \pi_1(N(x_o),x_o)$ we have that $\wh{\gamma}$
is a closed path, i.e., $\pi_1(N(\wh{x}_o),\wh{x}_o)\cong \pi_1(N(x_o),x_o)$.
\newtheorem{cor2}[prop]{Corollary}
\begin{cor2}
For all coverings $(\wh{M},q,M)$ with $Kerh_{\omega}\subseteq
\pi_1(\wh{M})$, we have that $N(\wh{x}_o)\cong N(x_o)$.
\end{cor2}
\pr These horizontal spaces are covered by that of the bundle induced over
the covering corresponding to $Kerh_{\omega}$. At the same time, they cover
$N(x_o)$. This, in view of the above proved proposition, completes the
proof. 
\begin{prop}
If we deal with $\wh{P}(\wh{M},G)$ where $\pi_1(\wh{M})=Kerh_{\omega}$, we
have $\pi_1(\wh{P})\cong p^{-1}_{\#}(Kerh_{\omega})$.
\end{prop}
\pr Since $(\wh{P},\wh{q},P)$ is a covering (cf. Lemma 1.1), it is enough
to see that $p^{-1}_{\#}(Kerh_{\omega})=\wh{q}_{\#}(\pi_1(\wh{P}))$. Let
$c$ be a path of $P$ closed at $x_o$, such that $[c]\in
p^{-1}_{\#}(Kerh_{\omega})$. As $[p\circ c]\in Kerh_{\omega}$, the path
$p\circ c$ is lifted to a closed path $\wh{p\circ c}$ on $\wh{M}$ (note
that $\pi_1(\wh{M})=Kerh_{\omega})$. Because $\wh{P}\cong G(\wh{x}_o)\times
\wh{M}$, the path $\wh{p\circ c}$ also gives the lifting of $c$ with
respect to $\wh{p}$. This implies that the lifting of $c$ with respect to
$\wh{q}$ is closed, hence $p^{-1}_{\#}(Kerh_{\omega})\subseteq
\wh{q}_{\#}(\pi_1(\wh{P}))$. As the other inclusion is obvious, the proof
is completed.

\smallskip
\indent
We finally identify $N(x_o)$ in terms of the holonomy morphism and data
related to $P(M,G)$.
\begin{prop}
Under the above notation $\pi_1(N(x_o),x_o)\cong
p^{-1}_{\#}(Kerh_{\omega})/$ $\pi_1(G (x_o),x_o)$.
\end{prop}
\pr Let $\gamma$ be a path of $P$, closed at $x_o$, such that
$p_{\#}([\gamma]) \in Kerh_{\omega}$. The horizontal lift $\ol{p\circ
\gamma}$ of this path with respect to $p$ is again closed at $x_o$ and
$[\ol{p \circ \gamma}] \in \pi_1 (N(x_o),x_o)$. As $p|N(x_o)$ is a covering
map, $\ol{p \circ \gamma}$ is nullhomotopic if and only if $[\gamma] \in
Kerp_{\#,x_o} \cong \pi_1(G(x_o),x_o)$ (cf. the exact homotopy sequence of
the fibration $G\rightarrow P \rightarrow M$). Hence the correspondence
$[\gamma] \rightarrow [\ol{p \circ \gamma}]$ defines a homomorphism of
$p^{-1}_{\#}(Kerh_{\omega})$ onto $\pi_1(N(x_o),x_o)$, whose kernel is
exactly the group $\pi_1(G(x_o),x_o)$.

\smallskip
\noindent
{\bf Added in Proof.} In a forthcoming work we consider similar problems in
the setting of associated fiber bundles and suitable flat connections
defined on them.

\vspace{1cm}

\noindent
{\sc University of Athens}

\noindent
{\sc Department of  Mathematics}

\noindent
{\sc Panepistimiopolis}

\noindent
{\sc GR 157 84, Athens}

\vspace{6pt}
\noindent
{\em E-mail address\/}: {\tt dlappas@cc.uoa.gr}

\end{document}